\documentclass[smallextended]{svjour3}
\usepackage{amssymb,amsmath}

\renewcommand*{\leq}{\leqslant} \renewcommand*{\geq}{\geqslant}

\newcommand*{\vare}{\varepsilon} \newcommand*{\varp}{\varphi}

\newcommand*{\bomega}{\text{\boldmath $\omega^2$}}

\newcommand*{\codim}{\mathop{\mathrm{codim}}\nolimits}
\newcommand*{\Fix}{\mathop{\mathrm{Fix}}\nolimits}

\newcommand*{\mQ}{\mathbb Q} \newcommand*{\mR}{\mathbb R}
\newcommand*{\mT}{\mathbb T} \newcommand*{\mZ}{\mathbb Z}

\newcommand*{\cG}{\mathcal G} \newcommand*{\cK}{\mathcal K}
\newcommand*{\cM}{\mathcal M} \newcommand*{\cP}{\mathcal P}
\newcommand*{\cT}{\mathcal T}

\newcommand*{\fG}{\mathfrak G}

\newcommand*{\tp}{\tilde{p}} \newcommand*{\tq}{\tilde{q}}
\newcommand*{\tu}{\tilde{u}} \newcommand*{\tv}{\tilde{v}}
\newcommand*{\tx}{\tilde{x}} \newcommand*{\ty}{\tilde{y}}
\newcommand*{\tz}{\tilde{z}}

\newcommand*{\hH}{\widehat{H}}

\newcommand*{\hcK}{\widehat{\cK}} \newcommand*{\hcM}{\widehat{\cM}}

\newcommand*{\DP}{\dot{p}} \newcommand*{\DQ}{\dot{q}}
\newcommand*{\DU}{\dot{u}} \newcommand*{\DV}{\dot{v}}
\newcommand*{\DX}{\dot{x}} \newcommand*{\DY}{\dot{y}}

\newcommand*{\DVARP}{\dot{\varp}}

\allowdisplaybreaks

\journalname{Arnold Mathematical Journal}

\begin{document}

\title{Integrable Hamiltonian Systems with a Periodic Orbit or Invariant Torus Unique in the Whole Phase Space}

\titlerunning{Integrable Hamiltonian Systems with a Unique Periodic Orbit}

\author{Mikhail B.~Sevryuk}

\institute{Mikhail B.~Sevryuk \at
V.~L.~Tal$'$roze Institute of Energy Problems of Chemical Physics, the Russian Academy of Sciences \\
Tel.: +7-499-1374104 \\
Fax.: +7-499-1378258 \\
\email{2421584@mail.ru, sevryuk@mccme.ru}
}

\date{Received: date / Accepted: date}

\maketitle

\begin{abstract}
It is very well known that periodic orbits of autonomous Hamiltonian systems are generically organized into smooth one-parameter families (the parameter being just the energy value). We present a simple example of an integrable Hamiltonian system (with an arbitrary number of degrees of freedom greater than one) with a unique periodic orbit in the phase space (which is not compact). Similar examples are given for Hamiltonian systems with a unique invariant torus (of any prescribed dimension) carrying conditionally periodic motions. Parallel examples for Hamiltonian systems with a compact phase space and with uniqueness replaced by isolatedness are also constructed. Finally, reversible analogues of all the examples are described.
\end{abstract}

\section{Introduction}
\label{Introduction}

Equilibrium points of autonomous Hamiltonian systems are generically isolated in the phase space, like equilibria of general dynamical systems. On the other hand, while generic periodic orbits of systems with dissipation are also isolated (and are therefore called limit cycles), periodic orbits of Hamiltonian systems are generically organized into smooth one-parameter families (the parameter being just the energy value). The proof is very simple \cite{AM78,HZ94}. Let $\gamma$ be a periodic orbit of an autonomous Hamiltonian system with $N$ degrees of freedom and $\fG$ the energy hypersurface containing $\gamma$. Consider a $(2N-2)$-dimensional local transversal section $\Sigma\subset\fG$ to $\gamma$ ($\Sigma\cap\gamma=\{O\}$) and the corresponding Poincar\'e map $\cP:\Sigma\to\Sigma$. Generically none of the eigenvalues of the linearization of $\cP$ at the fixed point $O$ is equal to one. If this is the case then, according to the implicit function theorem, any energy hypersurface sufficiently close to $\fG$ admits one and only one periodic orbit close to $\gamma$. This periodic orbit depends smoothly on the energy value.

By the way, this proof shows that periodic orbits of Hamiltonian systems with one degree of freedom are always included in smooth one-parameter families (each periodic orbit being a connected component of an energy line). Recall that the linearization of $\cP$ at $O$ is called the \emph{monodromy operator} of $\gamma$ (within $\fG$) and its eigenvalues (independent of $\Sigma$) are called the (characteristic) multipliers of $\gamma$. They occur in pairs $(\lambda,\lambda^{-1})$ of the same Jordan structure (and, in particular, of the same multiplicity).

The natural higher-dimensional generalization of the concept of a periodic orbit is an invariant torus carrying conditionally periodic motions. This is an invariant manifold $\cT$ diffeomorphic to the $n$-torus $\mT^n=\mR^n/2\pi\mZ^n$ and such that the dynamics on $\cT$ in a suitable angular coordinate $\varp\in\mT^n$ has the form $\DVARP=\omega$ with a constant vector $\omega\in\mR^n$ (called the frequency vector). Periodic orbits correspond to the case where $n=1$ and $\omega\neq 0$. If the frequencies $\omega_1,\ldots,\omega_n$ are incommensurable (rationally independent), one speaks of quasi-periodic motions. According to the Kolmogorov--Arnol$'$d--Moser (KAM) theory (see e.g.\ \cite{A89,AKN06,BHS96,KP03} and references therein), isotropic invariant tori of dimensions $2,\ldots,N$ carrying quasi-periodic motions are as typical for Hamiltonian systems with $N$ degrees of freedom as periodic orbits (recall that a submanifold of a symplectic manifold is said to be \emph{isotropic} if the restriction of the symplectic structure to this submanifold vanishes). Isotropic invariant $n$-tori ($2\leq n\leq N$) carrying quasi-periodic motions are generically organized into $n$-parameter families but these families are Cantor-like rather than smooth (however, they are foliated into smooth one-parameter subfamilies). In fact, the frequencies of invariant tori in the KAM theory are not merely incommensurable but strongly incommensurable (e.g., Diophantine), i.e., badly approximable by sets of commensurable frequencies.

The flow on $\mT^n$ afforded by the equation $\DVARP=\omega\in\mR^n$ is also said to be linear, parallel, rotational, or Kronecker. Therefore, invariant tori carrying conditionally periodic motions are sometimes called \emph{Kronecker tori} \cite{KP03}. A linear flow $g^t$ on $\mT^n$ with any frequency vector $\omega$ possesses the following \emph{recurrence property}: for any $T>0$ and $\vare>0$ there is $\tau>T$ such that for any $\varp\in\mT^n$ the distance between $\varp$ and $g^t(\varp)=\varp+t\omega$ is smaller than $\vare$. The distance here is to be understood with respect to, e.g., the flat Riemannian metric inherited from $\mR^n$. It is often convenient to extend the concept of a Kronecker torus to equilibria (invariant $0$-tori).

It is very well known that in various degenerate settings, periodic orbits as well as higher-dimensional isotropic Kronecker tori of a Hamiltonian system can constitute a family whose number of parameters is \emph{larger} than that in the generic case. Such situations are typical for superintegrable systems for which the number of independent first integrals exceeds the number of degrees of freedom (of course, not all of those integrals are pairwise in involution). For instance, consider the motion in the central force field in $\mR^N$ with a potential $V$ \cite{A89,AKN06}. This is a Hamiltonian system with $N$ degrees of freedom. If $V(r)=-c/r$ ($c>0$) then each trajectory with negative energy and non-zero angular momentum is an ellipse with a focus at the center of attraction (a Kepler ellipse). If $V(r)=cr^2$ ($c>0$) then each trajectory with non-zero angular momentum is an ellipse with the center at the center of attraction (a Hooke ellipse). In both the cases, an open domain of the phase space is foliated into periodic orbits.

Now the following question arises: can the number of parameters of a degenerate family of periodic orbits or higher-dimensional isotropic Kronecker tori of a Hamiltonian system be \emph{smaller} than that in the generic case? Can, for example, a periodic orbit $\gamma$ of a Hamiltonian system be isolated in the phase space? In other words, is it possible that there is a neighborhood $U$ of $\gamma$ such that $\gamma$ is the only periodic orbit entirely contained in $U$? Can a periodic orbit of a Hamiltonian system be unique in the whole phase space? There exists an extensive bibliography on periodic orbits of Hamiltonian systems (see e.g.\ \cite{HZ94} and references therein) but it seems that Hamiltonian systems with a unique periodic orbit have not been studied yet (on the contrary, Hamiltonian systems with no periodic orbits at all are a very popular subject of research). If $\gamma$ is the only periodic orbit of a Hamiltonian system then we arrive at the following astonishing picture: the energy hypersurface on which $\gamma$ lies contains only one periodic orbit (this is quite an ordinary situation, of course) but all the other energy hypersurfaces contain no periodic orbits!

In December 2017 and January 2018, the author and the user Khanickus of MathOverflow \cite{K18} constructed independently two very similar explicit (and exceedingly simple) examples of Hamiltonian systems in $\mR^4$ with a periodic orbit unique in the whole phase space. The main purpose of this short note is to present a generalization to the case where the dimension $n\geq 1$ of the invariant torus and the number $N\geq n+1$ of degrees of freedom are arbitrary.

\section{Invariant Tori in Hamiltonian Systems}
\label{Hamiltonian}

Let $n\geq 1$ and $m\geq 0$ be arbitrary integers and $\omega\in\mR^n$ an arbitrary vector. Our aim is to construct a Hamiltonian system on $\cM=\mR^{n+2m+2}\times\mT^n$ with a unique Kronecker $n$-torus, the frequency vector of this torus being $\omega$. Let $(u_1,\ldots,u_n,x,y,p_1,\ldots,p_m,q_1,\ldots,q_m)$ be coordinates in $\mR^{n+2m+2}$ and $(\varp_1,\ldots,\varp_n)$ angular coordinates in $\mT^n$. Consider the symplectic structure
\begin{equation}
\textstyle \bomega=\sum_{i=1}^ndu_i\wedge d\varp_i+dx\wedge dy+\sum_{j=1}^mdp_j\wedge dq_j
\label{structure}
\end{equation}
on $\cM$. The Hamilton function
\[
\textstyle H=\sum_{i=1}^n(\omega_iu_i+xu_i^2)+x^3/3+xy^2+\sum_{j=1}^m(p_j^3/3+p_jq_j^2)
\]
on the symplectic manifold $(\cM,\bomega)$ affords the equations of motion
\begin{equation}
\begin{aligned}
\DVARP_i &= \partial H/\partial u_i=\omega_i+2xu_i & &(1\leq i\leq n), \\
\DU_i &= -\partial H/\partial\varp_i=0 & &(1\leq i\leq n), \\
\DY &= \partial H/\partial x={\textstyle \sum_{i=1}^nu_i^2+x^2+y^2}, \\
\DX &= -\partial H/\partial y=-2xy, \\
\DQ_j &= \partial H/\partial p_j=p_j^2+q_j^2 & &(1\leq j\leq m), \\
\DP_j &= -\partial H/\partial q_j=-2p_jq_j & &(1\leq j\leq m).
\end{aligned}
\label{unique}
\end{equation}

The recurrence property of linear flows on tori and the equations $\DY=\sum_{i=1}^nu_i^2+x^2+y^2$ and $\DQ_j=p_j^2+q_j^2$ for $1\leq j\leq m$ imply that if a point $(u,x,y,p,q,\varp)\in\cM$ belongs to a Kronecker torus of~\eqref{unique} then $u=0$, $x=y=0$, and $p=q=0$. The $n$-torus
\begin{equation}
\cT=\{u=0, \;\; x=y=0, \;\; p=q=0\}
\label{torus}
\end{equation}
is indeed invariant under the flow of the Hamiltonian system~\eqref{unique} with $N=n+m+1$ degrees of freedom and carries conditionally periodic motions with the frequency vector $\omega$. It is therefore the only Kronecker $n$-torus of this system. Moreover, the torus~\eqref{torus} is isotropic and lies in the energy hypersurface $H^{-1}(0)$. If the frequencies $\omega_1,\ldots,\omega_n$ are incommensurable then the system~\eqref{unique} admits no other Kronecker tori whatsoever. If the rank of the set $\omega_1,\ldots,\omega_n$ over $\mQ$ is equal to $r<n$ then the system~\eqref{unique} has also Kronecker tori of dimensions $r,r+1,\ldots,n-1$ but they are contained in $\cT$.

For $n=1$ and $\omega\neq 0$ we obtain a Hamiltonian system with $m+2$ degrees of freedom and with a periodic orbit unique in the whole phase space. It is easy to see that the monodromy operator of $\cT$ within $H^{-1}(0)$ in this case is the identity operator in $\mR^{2m+2}$.

As was pointed out in the Introduction, there is no Hamiltonian system with one degree of freedom and with a unique periodic orbit in the phase space. The author does not know whether there exist Hamiltonian systems with $n$ degrees of freedom and with a unique (or just isolated) isotropic Kronecker $n$-torus in the phase space for $n\geq 2$.

The Hamiltonian system~\eqref{unique} is integrable: it admits $n+m+1$ first integrals
\begin{equation}
H, \quad u_i \;\; (1\leq i\leq n), \quad p_j^3/3+p_jq_j^2 \;\; (1\leq j\leq m)
\label{integrals}
\end{equation}
which are pairwise in involution and are functionally independent almost everywhere. However, on the torus~\eqref{torus} all these integrals vanish, and one has the degeneracy relations
\[
\textstyle dH=\sum_{i=1}^n\omega_idu_i, \quad d(p_j^3/3+p_jq_j^2)=0 \;\; (1\leq j\leq m).
\]

The common level surfaces of the integrals~\eqref{integrals} are not compact (each common level surface contains points with any value of $y\in\mR$), and this seems to be essential for the uniqueness of the torus. However, if one is interested just in the isolatedness, it is not hard to construct similar examples with compact common level surfaces of the integrals and even with a compact phase space. In fact, it suffices to take all the variables $(u_1,\ldots,u_n,x,y,p_1,\ldots,p_m,q_1,\ldots,q_m)$ modulo $2\pi$ in the example above. To be more precise, consider the symplectic manifold $\Bigl( \hcM,\bomega \Bigr)$ where $\hcM=\mT^{2n+2m+2}$ with angular coordinates
\[
(u_1,\ldots,u_n,\varp_1,\ldots,\varp_n,x,y,p_1,\ldots,p_m,q_1,\ldots,q_m)
\]
and the symplectic structure $\bomega$ is still given by the formula~\eqref{structure}. For any angular variable $z$ introduce the notation $\tz=\sin z$. The Hamilton function
\[
\textstyle \hH=\sum_{i=1}^n(\omega_i\tu_i+\tx\tu_i^2)+\tx^3/3+\tx\ty^2+\sum_{j=1}^m(\tp_j^3/3+\tp_j\tq_j^2)
\]
affords the equations of motion
\begin{equation}
\begin{aligned}
\DVARP_i &= \partial\hH/\partial u_i=\omega_i\cos u_i+\tx\sin 2u_i & &(1\leq i\leq n), \\
\DU_i &= -\partial\hH/\partial\varp_i=0 & &(1\leq i\leq n), \\
\DY &= \partial\hH/\partial x={\textstyle \bigl( \sum_{i=1}^n\tu_i^2+\tx^2+\ty^2 \bigr)\cos x}, \\
\DX &= -\partial\hH/\partial y=-\tx\sin 2y, \\
\DQ_j &= \partial\hH/\partial p_j=(\tp_j^2+\tq_j^2)\cos p_j & &(1\leq j\leq m), \\
\DP_j &= -\partial\hH/\partial q_j=-\tp_j\sin 2q_j & &(1\leq j\leq m).
\end{aligned}
\label{isolated}
\end{equation}

The $n$-torus~\eqref{torus} is again an isotropic Kronecker torus of the system~\eqref{isolated} with the frequency vector $\omega$. The torus $\cT$ lies in the energy hypersurface $\hH^{-1}(0)$. It is by no means unique; for instance, all the $2^{n+2m+2}$ tori given by the equations
\begin{equation}
\begin{gathered}
u_i=\delta_i \;\; (1\leq i\leq n), \\
x=\delta_{n+1}, \;\; y=\delta_{n+2}, \\
p_j=\delta_{n+2+j}, \;\; q_j=\delta_{n+m+2+j} \;\; (1\leq j\leq m)
\end{gathered}
\label{delta}
\end{equation}
where each of the numbers $\delta_1,\ldots,\delta_{n+2m+2}$ is equal to $0$ or to $\pi$ are isotropic Kronecker $n$-tori with the frequency vector $\omega$. Nevertheless, the torus $\cT$ is isolated in the phase space. Indeed, suppose that a point $(u,\varp,x,y,p,q)\in\hcM$ belongs to a Kronecker torus of~\eqref{isolated} entirely contained in the domain
\begin{equation}
\begin{gathered}
-\pi<u_i<\pi \;\; (1\leq i\leq n), \\
-\pi/2<x<\pi/2, \;\; -\pi<y<\pi, \\
-\pi/2<p_j<\pi/2, \;\; -\pi<q_j<\pi \;\; (1\leq j\leq m).
\end{gathered}
\label{pi}
\end{equation}
Then the recurrence property of linear flows on tori and the equations for $\DY$ and $\DQ_j$ ($1\leq j\leq m$) in~\eqref{isolated} imply that the point in question belongs to the torus~\eqref{torus}. Of course, all the inequalities in~\eqref{pi} are to be understood modulo $2\pi$: whenever $z$ is an angular variable, $a<z<b$ means that $z\in(a,b)\bmod 2\pi$.

Note that any neighborhood of the torus~\eqref{torus} has non-empty intersections with isotropic Kronecker tori of the system~\eqref{isolated} not contained in the domain~\eqref{pi}, for instance, with the $(n+1)$-tori $\{u=u^0, \;\; x=0, \;\; p=q=0\}$ where $u^0\neq 0$ is an arbitrary point sufficiently close to $0$. It is easy to verify that the frequency vector of such an $(n+1)$-torus is equal to
\[
\bigl( \omega_1\cos u^0_1,\ldots,\omega_n\cos u^0_n, \, [\zeta(\zeta+1)]^{1/2} \bigr), \quad {\textstyle \zeta=\sum_{i=1}^n\sin^2u^0_i>0}.
\]

The Hamiltonian system~\eqref{isolated} is also integrable: it admits $n+m+1$ first integrals
\[
\hH, \quad \tu_i \;\; (1\leq i\leq n), \quad \tp_j^3/3+\tp_j\tq_j^2 \;\; (1\leq j\leq m)
\]
which are pairwise in involution and are functionally independent almost everywhere. On the torus $\cT$ all these integrals vanish, and one has the degeneracy relations
\[
\textstyle d\hH=\sum_{i=1}^n\omega_id\tu_i, \quad d(\tp_j^3/3+\tp_j\tq_j^2)=0 \;\; (1\leq j\leq m).
\]

\section{Invariant Tori in Reversible Systems}
\label{reversible}

Many properties of Hamiltonian systems are also inherent (mutatis mutandis) in reversible systems. In particular, one can develop the reversible KAM theory which is in many respects parallel to the Hamiltonian KAM theory. The reader is referred to the books \cite{BHS96,S86} (and references therein) for the definition and main features of reversible dynamical systems. An invariant set of a reversible flow is usually said to be \emph{symmetric} if it is also invariant under the reversing involution. Consider an autonomous system reversible with respect to an involution $G$ such that the fixed point manifold $\Fix G$ of $G$ is not empty and all its connected components are of dimension $\ell$ and codimension $N$ (so that the phase space dimension is equal to $\ell+N$). Then symmetric equilibria of such a system are generically organized into smooth $(\ell-N)$-parameter families provided that $N\leq\ell$ (each family being just a smooth $(\ell-N)$-dimensional submanifold of $\Fix G$ consisting of equilibria), symmetric periodic orbits are generically organized into smooth $(\ell-N+1)$-parameter families provided that $N\leq\ell+1$, and symmetric invariant $n$-tori ($2\leq n\leq N$) carrying quasi-periodic motions with strongly incommensurable frequencies are generically organized into Cantor-like $(\ell-N+n)$-parameter families provided that $N\leq\ell+n-1$. Symmetric Kronecker $n$-tori with incommensurable frequencies and $n>N$ in such a system are impossible.

The Hamiltonian systems~\eqref{unique} and~\eqref{isolated} are reversible, the reversing involution in both the cases is given by the formula
\[
\cG:(u,\varp,x,y,p,q)\mapsto(u,-\varp,x,-y,p,-q),
\]
so that $\dim\Fix\cG=\codim\Fix\cG=n+m+1$. The Kronecker $n$-torus~\eqref{torus} in both the cases is symmetric.

These examples can be easily generalized to $G$-reversible systems with $\dim\Fix G\neq\codim\Fix G$ (the fact that such systems are natural to study was first emphasized by V.~I.~Arnol$'$d \cite{A84}). Let $n\geq 0$, $\ell\geq 0$, and $m\geq 0$ be arbitrary integers and $\omega\in\mR^n$ an arbitrary vector. Our aim is to construct a $G$-reversible system with $\dim\Fix G=\ell$, $\codim\Fix G=N=n+m+1$ and with a unique Kronecker $n$-torus, this torus being symmetric and the frequency vector of this torus being $\omega$. Note that if $n\leq 1$ then a $G$-reversible system with $\dim\Fix G>0$, $\codim\Fix G=n$ and with a unique symmetric Kronecker $n$-torus does definitely not exist. Indeed, if $\codim\Fix G=0$ and $\Fix G$ coincides with the phase space (i.e., if the involution $G$ is the identity transformation) then there is only one $G$-reversible vector field, namely, the zero field, and each point of the phase space is a symmetric equilibrium. On the other hand, it is easy to verify that if $\dim\Fix G=\ell$ and $\codim\Fix G=1$ then symmetric periodic orbits of any $G$-reversible system are always organized into smooth $\ell$-parameter families. The author does not know whether there are $G$-reversible systems with $\dim\Fix G>0$, $\codim\Fix G=n$ and with a unique (or just isolated) symmetric Kronecker $n$-torus for $n\geq 2$.

Let $(v_1,\ldots,v_\ell,y,q_1,\ldots,q_m)$ be coordinates in $\mR^{\ell+m+1}$ and $(\varp_1,\ldots,\varp_n)$ angular coordinates in $\mT^n$. The dynamical system
\begin{equation}
\begin{aligned}
\DVARP_i &= \omega_i & &(1\leq i\leq n), \\
\DV_k &= 0 & &(1\leq k\leq\ell), \\
\DY &= {\textstyle \sum_{k=1}^\ell v_k^2+y^2+\sum_{j=1}^mq_j^2}, \\
\DQ_j &= 0 & &(1\leq j\leq m)
\end{aligned}
\label{revunique}
\end{equation}
on $\cK=\mT^n\times\mR^{\ell+m+1}$ is reversible with respect to the phase space involution
\begin{equation}
G:(\varp,v,y,q)\mapsto(-\varp,v,-y,-q),
\label{involution}
\end{equation}
here $\dim\Fix G=\ell$ and $\codim\Fix G=n+m+1$.

The recurrence property of linear flows on tori and the equation for $\DY$ in~\eqref{revunique} imply that if a point $(\varp,v,y,q)\in\cK$ belongs to a Kronecker torus of~\eqref{revunique} (symmetric or not) then $v=0$, $y=0$, and $q=0$. The $n$-torus
\begin{equation}
\{v=0, \;\; y=0, \;\; q=0\}
\label{revtorus}
\end{equation}
is indeed Kronecker. It is therefore the only Kronecker $n$-torus of the system~\eqref{revunique}. Moreover, the torus~\eqref{revtorus} is symmetric, and its frequency vector is equal to $\omega$.

As in the Hamiltonian case, this example can be compactified. Let
\[
(\varp_1,\ldots,\varp_n,v_1,\ldots,v_\ell,y,q_1,\ldots,q_m)
\]
be angular coordinates in $\hcK=\mT^{n+\ell+m+1}$. The dynamical system
\begin{equation}
\begin{aligned}
\DVARP_i &= \omega_i & &(1\leq i\leq n), \\
\DV_k &= 0 & &(1\leq k\leq\ell), \\
\DY &= {\textstyle \sum_{k=1}^\ell\tv_k^2+\ty^2+\sum_{j=1}^m\tq_j^2}, \\
\DQ_j &= 0 & &(1\leq j\leq m)
\end{aligned}
\label{revisolated}
\end{equation}
on $\hcK$ (as before, here $\tz=\sin z$ for any angular variable $z$) is reversible with respect to the phase space involution $G$ given by the formula~\eqref{involution}, and $\dim\Fix G=\ell$, $\codim\Fix G=n+m+1$.

The $n$-torus~\eqref{revtorus} is again a symmetric Kronecker torus of the system~\eqref{revisolated} with the frequency vector $\omega$. It is not unique, of course; for instance, all the $2^{\ell+m+1}$ tori given by the equations
\[
v_k=\delta_k \;\; (1\leq k\leq\ell), \quad y=\delta_{\ell+1}, \quad q_j=\delta_{\ell+1+j} \;\; (1\leq j\leq m)
\]
where each of the numbers $\delta_1,\ldots,\delta_{\ell+m+1}$ is equal to $0$ or to $\pi$ (cf.~\eqref{delta}) are symmetric Kronecker $n$-tori with the frequency vector $\omega$. However, the torus~\eqref{revtorus} is isolated in the phase space. Indeed, suppose that a point $(\varp,v,y,q)\in\hcK$ belongs to a Kronecker torus of~\eqref{revisolated} (symmetric or not) entirely contained in the domain
\begin{equation}
-\pi<v_k<\pi \;\; (1\leq k\leq\ell), \quad -\pi<y<\pi, \quad -\pi<q_j<\pi \;\; (1\leq j\leq m)
\label{revpi}
\end{equation}
(cf.~\eqref{pi}). Then the recurrence property of linear flows on tori and the equation for $\DY$ in~\eqref{revisolated} imply that the point in question belongs to the torus~\eqref{revtorus}.

Note that any neighborhood of the torus~\eqref{revtorus} has non-empty intersections with symmetric Kronecker tori of the system~\eqref{revisolated} not contained in the domain~\eqref{revpi}, for instance, with the $(n+1)$-tori $\{v=v^0, \;\; q=0\}$ where $v^0\neq 0$ is an arbitrary point sufficiently close to $0$. It is easy to see that the frequency vector of such an $(n+1)$-torus is equal to
\[
\bigl( \omega_1,\ldots,\omega_n, \, [\xi(\xi+1)]^{1/2} \bigr), \quad {\textstyle \xi=\sum_{k=1}^\ell\sin^2v^0_k>0}.
\]

\begin{acknowledgements}
The author is grateful to L.~M.~Lerman for useful discussions on the role of compactness of various invariant manifolds in the phenomena studied in this note.
\end{acknowledgements}

\end{document}